\documentclass[dvips,aos]{imsart}

\usepackage{amsfonts}

\RequirePackage[OT1]{fontenc}
\RequirePackage{amsthm,amsmath,natbib}
\RequirePackage[colorlinks]{hyperref}
\RequirePackage{hypernat}

\startlocaldefs

\numberwithin{equation}{section}
\theoremstyle{plain}
\newtheorem{thm}{Theorem}[section]
\newtheorem{cor}[thm]{Corollary}
\newtheorem{lemma}[thm]{Lemma}
\newtheorem{proposition}[thm]{Proposition}

\def\triangleq{\stackrel{\triangle}{=}}

\def\argmin{\mathop{\mathrm{arg\,min}}}
\def\sf{}

\def\wideth{\widehat}
\def\scm{{\cal M}}

\def\wh{\widehat}
\def\eps{\varepsilon}
\def\A{\mbox{$\mathcal A$}}
\def\B{\mbox{$\mathcal B$}}

\def\Z{\mbox{$\mathcal Z$}}
\def\RR{\mathbb R}

\def\PP{\mathbb{P}}

\def\RE#1#2{$\text{RE}(#1,#2)$}
\def\REm#1#2#3{$\text{RE}(#1,#2,#3)$}

\newcommand{\bD}{\mbox{\boldmath$\delta$}}
\newcommand{\bb}{\mbox{\boldmath$\beta$}}
\def\bb{\beta}
\newcommand{\bmf}{\mbox{\boldmath$f$}}
\newcommand{\bmy}{\mbox{\boldmath$y$}}
\newcommand{\bmw}{\mbox{\boldmath$w$}}

\def\argmin{\mathop{\rm arg\, min}}

\endlocaldefs

\begin{document}

\begin{frontmatter}
\title{Simultaneous analysis of Lasso and Dantzig selector\protect\thanksref{T1}}
\runtitle{Lasso and Dantzig selector}
\thankstext{T1}{Partially supported by NSF grant DMS-0605236,
ISF Grant, France-Berkeley Fund, the grant ANR-06-BLAN-0194 and the
European Network of Excellence PASCAL. }


\begin{aug}
\author{\fnms{Peter J.} \snm{ Bickel}
\ead[label=e1]{bickel@stat.berkeley.edu}},
\author{\fnms{Ya'acov} \snm{Ritov}
\ead[label=e2]{yaacov.ritov@gmail.com}}
\and
\author{\fnms{Alexandre B.} \snm{Tsybakov}
\ead[label=e3]{tsybakov@ccr.jussieu.fr}
\ead[label=u3,url]{http://www.foo.com}}

\runauthor{Bickel et al.}

\affiliation{}

\address{Department of Statistics\\University of California at Berkeley,\\ CA USA\\
\printead{e1}
\phantom{E-mail:\ }}

\address{Jerusalem, Israel\\
\printead{e2}
\phantom{E-mail:\ }}

\address{Laboratoire de Probabilit\'es et
Mod\`eles Al\'eatoires\\Universit\'e Paris VI, France.\\
\printead{e3}\\
}
\end{aug}

\begin{abstract}
We show that, under a sparsity scenario, the Lasso estimator and the
Dantzig selector exhibit similar behavior. For both methods we
derive, in parallel, oracle inequalities for the prediction risk in
the general nonparametric regression model, as well as bounds on the
$\ell_p$ estimation loss for $1\le p\le 2$ in the linear model when
the number of variables can be much larger than the sample size.
\end{abstract}



\end{frontmatter}

\def\citenumfont{}

\section{Introduction}\label{sec0}

During the last few years a great deal of attention has been focused
on the $\ell_1$ penalized least squares (Lasso) estimator of
parameters in high-dimensional linear regression when the number of
variables can be much larger than the sample size
\cite{det04,ehjt04,kf00,mb06,my06,os00a,os00b,t96,zh06,zy05}. Quite
recently, Candes and Tao \cite{ct07} have proposed a new estimate
for such linear models, the Dantzig selector, for which they
establish optimal $\ell_2$ rate properties under a sparsity
scenario, i.e., when the number of non-zero components of the true
vector of parameters is small.

 Lasso estimators have been
also studied in the nonparametric regression setup
\cite{jn00,n00,gr04,btw04,btw06a,btw07,btw06b}. In particular, Bunea
et al. \cite{btw04,btw06a,btw07,btw06b} obtain sparsity oracle
inequalities for the prediction loss in this context and point out
the implications for minimax estimation in classical non-parametric
regression settings, as well as for the problem of aggregation of
estimators. An analog of Lasso for density estimation with similar
properties (SPADES) is proposed in \cite{btw07a}. Modified versions
of Lasso estimators  (non-quadratic   terms and/or
penalties slightly different from $\ell_1$) for nonparametric
regression with random design are suggested and studied under
prediction loss in \cite{k06,vdg06}. Sparsity oracle inequalities
for the Dantzig selector with random design are obtained in \cite{k07}.
In linear fixed design regression, Meinshausen and Yu \cite{my06}
establish a bound on the $\ell_2$ loss for the coefficients of Lasso
which is quite different from the bound on the same loss for the Dantzig selector proven in
\cite{ct07}.

The main message of this paper is that under a sparsity scenario,
the Lasso and the Dantzig selector exhibit similar behavior, both
for linear regression and for nonparametric regression models, for
$\ell_2$ prediction loss and for $\ell_p$ loss in the coefficients
for $1\le p\le 2$. All the results of the paper are non-asymptotic.

Let us specialize to the case of linear regression with many
covariates, $\bmy=X\bb+\bmw$ where $X$ is the $n\times M$
deterministic design matrix, with $M$ possibly much larger than $n$,
and $\bmw$ is a vector of i.i.d. standard normal random variables.
This is the situation considered most recently by Candes and Tao
\cite{ct07} and Meinshausen and Yu \cite{my06}.  Here sparsity
specifies that the high-dimensional vector $\bb$ has coefficients
that are mostly 0.

We develop general tools to study these two estimators in parallel.
For the fixed design Gaussian regression model we recover, as
particular cases, sparsity oracle inequalities for the Lasso, as in
Bunea et al. \cite{btw07}, and $\ell_2$ bounds for the coefficients
of Dantzig selector, as in Candes and Tao \cite{ct07}. This is
obtained as a consequence of our more general results:
\begin{itemize}
\item  In the nonparametric regression model, we prove sparsity oracle
    inequalities for the Dantzig
    selector, that is, bounds on the prediction loss in terms of
    the best possible (oracle) approximation under the sparsity
    constraint.
\item Similar sparsity oracle inequalities are proved for the Lasso in
    the nonparametric regression model, and this is done under more
    general assumptions on the design matrix than in \cite{btw07}.
\item We prove that, for nonparametric regression, the Lasso and  the
    Dantzig selector are approximately equivalent in terms of the
    prediction loss.
\item We develop geometrical assumptions which are considerably weaker
    than those of Candes and Tao \cite{ct07} for the Dantzig selector and
    Bunea et al.\cite{btw07} for the Lasso. In the context of linear
    regression where the number of variables is possibly much larger than
    the sample size these assumptions imply the result of \cite{ct07} for
    the $\ell_2$ loss and generalize it to $\ell_p$ loss, $1\leq p\leq2$,
    and to prediction loss. Our bounds for the Lasso differ from those
    for Dantzig selector only in numerical constants.
\end{itemize}
We begin, in the next section, by defining the Lasso and Dantzig
procedures and the notation.  In Section \ref{sec2} we present our
key geometric assumptions. Some sufficient conditions for these
assumptions are given in Section \ref{discRE}, where they are also
compared to those of \cite{ct07} and \cite{my06} as well as to ones
appearing in \cite{btw07} and \cite{btw06b}. We  note a weakness of
our assumptions, and hence of those in the papers we cited, and we
discuss a way of slightly remedying them. Sections \ref{sec3},
\ref{sec5} give some equivalence results and sparsity oracle
inequalities for the Lasso and Dantzig estimators in the general
nonparametric regression model. Section \ref{sec6} focuses on the
linear regression model and includes a final discussion. Two
important technical lemmas are given in Appendix \ref{app:proofs} as
well as most of the proofs.

\section{Definitions and notation}\label{sec1}

Let $ (Z_1,Y_1),$ $\ldots,(Z_n,Y_n) $ be a sample of independent random
pairs with
\begin{equation*}
Y_i = f(Z_i) + W_i, \quad i=1,\ldots,n,
\end{equation*}
where $f:{\mathcal Z} \to \RR$ is an unknown regression function to
be estimated, ${\mathcal Z}$ is a Borel subset of $\RR^d$, the
$Z_i$'s are fixed elements in $\Z$ and the regression errors $W_i$
are Gaussian. Let $\mathcal{F}_M = \{f_1, \ldots, f_M\}$ be a finite
dictionary of functions $f_j:{\mathcal Z} \to \RR$, $j=1,\dots,M$.
We assume throughout that $M\ge2$.

Depending on the statistical targets, the dictionary $\mathcal{F}_M$
can contain qualitatively different parts. For instance, it can be a
collection of basis functions used to approximate $f$ in the
nonparametric regression model (e.g., wavelets, splines with fixed
knots, step functions). Another example is related to the
aggregation problem where the $f_j$ are estimators arising from $M$
different methods. They  can also correspond  to $M$ different
values of the tuning parameter of the same method. Without much loss
of generality, these estimators $f_j$ are treated as fixed
functions: the results are viewed as being conditioned on the sample
the $f_j$ are based on.

The selection of the dictionary can be very important to make the
estimation of $f$ possible.  We assume implicitly that $f$ can be
well approximated by a member of the span of $\mathcal{F}_M$.
However this is not enough.  In this paper, we have in mind the
situation where $M\gg n$, and $f$ can be estimated reasonably only
because it can approximated by a linear combination of a small
number of members of $\mathcal{F}_M$, or in other words, it has a
sparse approximation in the span of ${\mathcal{F}}_M$. But when
sparsity is an issue, equivalent bases can have different
properties: a function which has a sparse  representation in one
basis may not have it in another one, even if both of them span the
same linear space.

Consider the matrix $X= (f_j(Z_i))_{i,j}$, $i=1,\dots,n$,
$j=1,\dots,M$ and the vectors  $\bmy=(Y_1,\dots,Y_n)^T$, ${\bmf
}=(f(Z_1),\dots,f(Z_n))^T$, $\bmw=(W_1,\dots,W_n)^T$.  With this
notation,
\begin{eqnarray*}
\bmy &= {\bmf} +\bmw.
\end{eqnarray*} We will write
$|x|_p$ for the $\ell_p$ norm of $x\in\RR^M$, $1\le p\le\infty$. The notation  $\|\cdot\|_n$ stands for the empirical norm:
$$\|g\|_n= \sqrt{\frac1n \sum_{i=1}^n g^2(Z_i)}$$
for any $g:\mathcal Z\to\RR$. We suppose that $\|f_j\|_n\neq 0$, $j=1,\dots, M$.
Set
$$
f_{\max} = \max_{1\le j\le M} \|f_j\|_n\,, \quad f_{\min}= \min_{1\le j\le M} \|f_j\|_n\,.
$$

 For any
$\beta= (\beta_1,\dots,\beta_M)\in \RR^M$, define ${\sf
f}_{\beta} = \sum_{j=1}^M \beta_j f_j$, or explicitly, ${\sf
f}_{\beta}(z) = \sum_{j=1}^M \beta_j f_j(z)$, and $\bmf_\beta =X\beta$.
The estimates we consider are all of the form ${ f}_{\tilde \beta} (\cdot)$
where $\tilde \beta$ is data determined.
Since we consider mainly sparse vectors $\tilde \beta$, it will be convenient to define the following. Let $$ \scm(\beta)= \sum_{j=1}^M I_{\{ \beta_j\ne 0\}} =
|J(\beta)|$$ denote the number of non-zero coordinates of
$\beta$, where $I_{\{ \cdot\}}$ denotes the indicator function,
$J(\beta)=\left \{j\in\{1,\ldots,M\} \right.:$ $\left. \beta_j
\ne 0\right\}$, and $|J|$ denotes the cardinality of $J$. The value
$\scm(\beta)$ characterizes the {\it sparsity} of the vector
$\beta$: the smaller $\scm(\beta)$, the ``sparser"
 $\beta$.
For a vector $\bD\in \RR^M$ and a subset $J\subset \{1,\dots,M\}$ we
denote by $\bD_J$ the vector in $\RR^M$ which has the same
coordinates as $\bD$ on $J$ and zero coordinates on the complement
$J^c$ of $J$.

 Introduce the residual sum of squares
$$ \wh{S} (\beta) = \frac1n \sum_{i=1}^n \{ Y_i - {\sf
f}_{\beta}(Z_i) \}^2,$$ for all $\beta\in\RR^M$.  Define the
Lasso solution $\wh{\beta}_L =
(\wh{\beta}_{1,L},\dots,\wh{\beta}_{M,L})$ by
\begin{equation}\begin{split}\label{lass}
\wh{\beta}_L = \argmin_{\beta\in\RR^M} \left\{ \wh{S}(\beta) +
2r\sum_{j=1}^M  \|f_j\|_n |\beta_j| \right\},
\end{split}\end{equation}
where $r>0$ is some tuning constant, and introduce the corresponding Lasso
estimator
\begin{equation}\begin{split}
\label{lse} \widehat{f}_L(x)={\sf
f}_{\wh{\beta}_L}(x) =\sum_{j=1}^M \wh{\beta}_{j,L} f_j(z).
\end{split}\end{equation}

The criterion in (\ref{lass}) is convex in $\beta$, so that standard
convex optimization procedures can be used to compute
$\wh{\beta}_L$. We refer to
\cite{os00a,os00b,ehjt04,t05,fjht07,mvb08} for detailed discussion
of these optimization problems and fast algorithms.

A necessary and sufficient condition of the minimizer in
(\ref{lass}) is that 0 belongs to the subdifferential of the convex
function $\beta\mapsto n^{-1}|y-X\beta|_2^2 + 2r|D^{1/2}\beta|_1$.
This implies that the Lasso selector $\widehat \beta_L$ satisfies
the constraint:
\begin{equation}\begin{split}\label{dc}
\Big|\frac{1}{n}D^{-1/2}X^T(y-X\widehat\beta_L)\Big|_\infty \le r,
\end{split}\end{equation}
where $D$ is the diagonal matrix
$$
D = {\rm diag}\{\|f_1\|_n^{2},\dots,\|f_M\|_n^{2}\}.
$$
More generally, we will say that $\beta\in\RR^M$ satisfies the
Dantzig constraint if $\beta$ belongs to the set
$$
\left\{\beta\in\RR^M:\,
\Big|\frac{1}{n}D^{-1/2}X^T(y-X\beta)\Big|_\infty \le r \right\}\,.
$$

The Dantzig estimator of the regression function $f$ is based on a
particular solution of \eqref{dc}, the Dantzig selector
$\wh{\beta}_D$ which is defined as a vector having the smallest
$\ell_1$ norm among all $\beta$ satisfying the Dantzig constraint:
\begin{equation}\begin{split}\label{da2}
\wh{\beta}_D = \argmin\Big\{|\beta|_1: \
\Big|\frac{1}{n}D^{-1/2}X^T(y-X\beta)\Big|_\infty \le r \Big\}.
\end{split}\end{equation}
The Dantzig estimator is defined by
\begin{equation}\begin{split}\label{da1}
\widehat{f}_D(z)={\sf f}_{\wh{\beta}_D}(z) =\sum_{j=1}^M
\wh{\beta}_{j,D} f_j(z).
\end{split}\end{equation}
where $\wh{\beta}_D = (\wh{\beta}_{1,D},\dots,\wh{\beta}_{M,D})$ is
the Dantzig selector. By the definition of Dantzig selector, we have
$|\widehat\beta_D|_1\le |\widehat\beta_L|_1$.

The Dantzig selector is computationally feasible, since it reduces to a
linear programming problem \cite{ct07}.

Finally for any $n\ge 1$, $M\ge2$, we consider the Gram matrix $$ \Psi_{n} =
\frac1n X^TX = \left( \frac1n \sum_{i=1}^n f_j(Z_i) f_{j^\prime}
(Z_i) \right)_{ 1\le j,j^\prime \le M},$$
and let $\phi_{\max}$ denote  the maximal eigenvalue of $\Psi_n$.

\section{Restricted eigenvalue assumptions}\label{sec2}

We now introduce the key assumptions on the Gram matrix that are
needed to guarantee nice statistical properties of Lasso and Dantzig
selector. Under the sparsity scenario we are typically interested in
the case where $M>n$, and even $M\gg n$. Then the matrix $\Psi_n$ is
degenerate, which can be written as
$$
\min_{\bD\in \RR^M: \bD \neq 0} \ \
\frac{(\bD^T\Psi_{n}\bD)^{1/2}}{|\bD|_2} \, \equiv
 \min_{\bD\in \RR^M: \bD \neq 0}
\ \ \frac{|X\bD|_2}{\sqrt{n}|\bD|_2} \,= \, 0.
$$
Clearly,  ordinary least squares does not work in this case,
since it requires positive definiteness of $\Psi_n$, i.e.
\begin{equation}\label{31}
\min_{\bD\in \RR^M: \bD \neq 0}
\ \ \frac{|X\bD|_2}{\sqrt{n}|\bD|_2} \,> \, 0.
\end{equation}
It turns out that the Lasso and Dantzig selector require much weaker
assumptions: the minimum  in (\ref{31}) can be replaced by
the minimum over a restricted set of vectors, and the norm  $|\bD|_2$ in
the denominator of the condition can be replaced by the $\ell_2$ norm of only a part of $\bD$.

One of the properties of both the Lasso and the Dantzig selector is
that, for the linear regression model, the residuals $\bD=\hat
\beta_L-\beta$ and $\bD=\hat \beta_D-\beta$ satisfy, with
probability close to 1,
\begin{equation}
\label{c0ins} |\bD_{J_0^c}|_1\leq c_0|\bD_{J_0}|_1
\end{equation}
where $J_0=J(\beta)$ is the set of non-zero coefficients of the true
parameter $\beta$ of the model. For the linear regression model, the
vector of Dantzig residuals $\bD$ satisfies \eqref{c0ins} with
probability close 1 if $c_0=1$ and $M$ is large (cf. (\ref{ll22})
and the fact that $\beta$ of the model satisfies the Dantzig
constraint with probability close 1 if $M$ is large). A similar
inequality holds for the vector of Lasso residuals $\bD=\wh \beta_L
- \beta$, but this time with $c_0=3$, cf. Corollary \ref{corx}.

Now, consider for example, the case where the elements of the Gram
matrix $\Psi_n$ are close to those of a positive definite $M\times
M$ matrix $\Psi$. Denote by
$\eps_n\triangleq\max_{i,j}|(\Psi_{n}-\Psi)_{i,j}|$ the maximal
difference between the elements of the two matrices. Then for any
$\bD $ satisfying \eqref{c0ins} we get
\begin{equation}
\label{REheur}
\begin{split}
\frac{\bD^T\Psi_n\bD}{|\bD|_2^2} &= \frac{\bD^T\Psi\bD
+\bD^T(\Psi_n-\Psi)\bD}{|\bD|_2^2}
\\
&\geq \frac{\bD^T\Psi\bD}{|\bD|_2^2} - \frac{\eps_n
|\bD|_1^2}{|\bD|_2^2}
\\
&\geq \frac{\bD^T\Psi\bD}{|\bD|_2^2} - \eps_n
\left(\frac{(1+c_0)|\bD_{J_0}|_1}{|\bD_{J_0}|_2} \right)^2
\\
&\geq \frac{\bD^T\Psi\bD}{|\bD|_2^2} - \eps_n (1+c_0)^2|J_0|.
\end{split}
\end{equation}
Thus, for $\bD $ satisfying \eqref{c0ins} which are the vectors that
we have in mind, and for $\eps_n|J_0|$ small enough, the LHS of
\eqref{REheur} is bounded away from 0. It means that we have a kind
of ``restricted" positive definiteness which is valid only for the
vectors satisfying \eqref{c0ins}. This suggests the following
conditions that will suffice for the main argument of the paper. We
refer to these conditions as {\it restricted eigenvalue} (RE)
assumptions.

Our first RE assumption is:

\begin{description}
\item[Assumption \RE s{c_0}:]
For some integer $s$ such that $1\le s\le M$,  and a positive number
$c_0$ the following condition holds:
$$
\kappa(s, c_0) \triangleq \min_{\substack{J_0\subseteq \{1,\dots,M\},\\
\\ |J_0|\le s}} \ \ \min_{\substack{\bD\neq 0, \\ |\bD_{J_0^c}|_1\le c_0|\bD_{J_0}|_1}}
\ \ \frac{|X\bD|_2}{\sqrt{n}|\bD_{J_0}|_2} \,> \, 0.
$$

\end{description}

The integer $s$ here plays the role of an upper bound on the
sparsity $\scm(\beta)$ of a vector of coefficients $\beta$.

Note that if Assumption \RE s{c_0} is satisfied with $c_0\geq 1$,
then $$\min\{|X\bD|_2:\; \scm(\bD)\le 2s,\bD\ne 0\}>0.$$ In words,
the square submatrices of size $\le 2s$ of the Gram matrix are
necessarily positive definite. Indeed, suppose that for some $\bD\ne
0$ we have simultaneously $\scm(\bD)\leq 2s$ and $X\bD=0$. Partition
$J(\bD)$ in two sets: $J(\bD)=I_0\cup I_1$, such that $|I_i|\leq s$,
$i=0,1$. Without loss of generality, suppose that $|\bD_{I_1}|_1\le
|\bD_{I_0}|_1$. Since, clearly, $|\bD_{I_1}|_1=|\bD_{I_0^c}|_1$ and
$c_0\ge 1$, we have $|\bD_{I_0^c}|_1\le c_0|\bD_{I_0}|_1$. Hence
$\kappa(s,c_0)=0$, a contradiction.

To introduce the second assumption we need some more notation. For
integers $s,m$ such that $1\le s \le M/2$ and $m\ge s$, $s+m\le M$,
a vector $\bD\in \RR^M$ and a set of indices $J_0\subseteq
\{1,\dots,M\}$ with $|J_0|\le s$, denote by $J_{1}$ the subset of
$\{1,\dots,M\}$ corresponding to the $m$ largest in absolute value
coordinates of $\bD$ outside of $J_0$ and define $J_{01}\triangleq
J_0\cup J_1$. Clearly, $J_{1}$ and $J_{01}$ depend on $m$ but we do
not indicate this in our notation for the sake of brevity.

\begin{description}
\item[Assumption \REm sm{c_0}:]
$$
\kappa(s,m,c_0) \triangleq \min_{\substack{J_0\subseteq \{1,\dots,M\},\\
\\ |J_0|\le s}} \ \ \min_{\substack{\bD\neq 0,\\ |\bD_{J_0^c}|_1\le c_0|\bD_{J_0}|_1}}
\ \ \frac{|X\bD|_2}{\sqrt{n}|\bD_{J_{01}}|_2} \,> \, 0.
$$
\end{description}

Note that the only difference between the two assumptions is in the
denominators, and $\kappa(s,m,c_0)\leq \kappa(s,c_0)$. As written,
for fixed $n$, the two assumptions are equivalent. However,
asymptotically for large $n$, Assumption \RE s{c_0} is less
restrictive than \REm sm{c_0}, since the ratio $\kappa(s,m,c_0)/
\kappa(s,c_0)$ may tend to 0 if $s$ and $m$ depend on $n$. For our
bounds on the prediction loss and on the $\ell_1$ loss of the Lasso
and Dantzig estimators we will only need Assumption \RE s{c_0}.
Assumption \REm sm{c_0} will be required exclusively for the bounds
on the $\ell_p$ loss with $1<p\le 2$.

Note also that Assumptions \RE {s'}{c_0} and \REm {s'}m{c_0}
imply Assumptions \RE s{c_0} and \REm sm{c_0}
respectively if $s'>s$.

\section{Discussion of the RE assumptions}
\label{discRE}

There exist several simple sufficient conditions for Assumptions \RE
s{c_0} and \REm sm{c_0} to hold. Here we discuss some of them.

For a real number $1\le u\le M$ we introduce the following
quantities that we will call {\it restricted eigenvalues}:
\begin{equation*}
\begin{split}
\phi_{\min}(u) &= \min_{x\in \RR^M: 1\le \scm(x) \le u }
\frac{x^T\Psi_n
x}{|x|_2^2}, \\
 \phi_{\max}(u) &= \max_{x\in \RR^M: 1\le \scm(x) \le
u } \frac{x^T\Psi_n x}{|x|_2^2}.
\end{split}
\end{equation*}
Denote by $X_J$ the $n\times |J|$ submatrix of $X$ obtained by
removing from $X$ the columns that do not correspond to the indices
in $J$, and for $1\le m_1,m_2\le M$ introduce the following
quantities called {\it restricted correlations}:
\begin{equation*}
\theta_{m_1,m_2} = \max\left\{ \frac{
c_1^TX_{I_1}^TX_{I_2}c_2}{n\,|c_1|_2\,|c_2|_2}: I_1\cap
I_2=\emptyset, |I_i|\leq m_i, c_i\in\RR^{I_i}\setminus \{0\}, i=1,2
\right\}.
\end{equation*}

In Lemma \ref{l2} below we show that a sufficient condition for \RE
s{c_0} and \REm ss{c_0} to hold is given, for example, by the
following assumption on the Gram matrix.

\begin{description}
\item[Assumption 1:]Assume
$$
\phi_{\min}(2s) > c_0\theta_{s,2s}
$$
for some integer $1\le s \le M/2$ and a constant $c_0>0$.
\end{description}

This condition with $c_0=1$ appeared in \cite{ct07}, in connection
with the Dantzig selector. Assumption 1 is more general: we can have
here an arbitrary constant $c_0>0$ which will allow us to cover not
only the Dantzig selector but also the Lasso estimators, and to
prove oracle inequalities for the prediction loss when the model is
nonparametric.

Our second sufficient condition for \RE s{c_0} and \REm sm{c_0} does
not need bounds on correlations. Only bounds on the minimal and
maximal eigenvalues of ``small" submatrices of the Gram matrix
$\Psi_n$ are involved.

\begin{description}
\item[Assumption 2:] Assume
$$
m\phi_{\min}(s+m) > c_0^2 s\phi_{\max}(m)
$$
for some integers $s,m$ such that $1\le s \le M/2$, $m\ge s$, and
$s+m\le M$, and a constat $c_0>0$.
\end{description}

Assumption 2 can be viewed as a weakening of the condition on
$\phi_{\min}$ in \cite{my06}. Indeed, taking $s+m=s\log n$ (we
assume w.l.o.g. that $s\log n$ is an integer and $n>3$) and assuming
that $\phi_{\max}(\cdot)$ is uniformly bounded by a constant we get
that Assumption 2 is equivalent to
\begin{equation*}
\phi_{\min}(s \log n) > c/\log n
\end{equation*}
where $c>0$ is a constant. The corresponding slightly stronger
assumption in \cite{my06} is stated in asymptotic form (for
$s=s_n\to\infty$):
$$
\liminf_{n}\phi_{\min}(s_n \log n) > 0.
$$

The following two constants are useful when Assumptions 1 and 2 are
considered:
$$
\kappa_1 (s,c_0) = \sqrt{\phi_{\min}(2s)}\left(1-
\frac{c_0\,\theta_{s,2s}}{\phi_{\min}(2s)}\ \right)\,
$$
and
$$
\kappa_2 (s,m,c_0) = \sqrt{\phi_{\min}(s+m)}\left(1-
c_0\sqrt{\frac{s\,\phi_{\max}(m)}{m\,\phi_{\min}(s+m)}}\ \right)\, .
$$
The next lemma shows that if Assumptions 1 or 2 are satisfied, then
the quadratic form $x^T\Psi_n x$ is positive definite  on some
restricted sets of vectors $x$. The construction of the lemma is
inspired by Candes and Tao \cite{ct07} and covers, in particular,
the corresponding result in \cite{ct07}.

\begin{lemma}\label{l2} Fix an integer $1\le s \le M/2$ and a
constant $c_0>0$.

(i) Let Assumption 1 be satisfied. Then Assumptions \RE s{c_0} and
\REm ss{c_0} hold with $\kappa (s,c_0) =\kappa (s,s,c_0)=\kappa_1
(s,c_0)$. Moreover, for any subset $J_0$ of $\{1,\dots,M\}$ with
cardinality $|J_0|\le s$, and any $\bD\in \RR^M$ such that
\begin{equation}\label{l21}
|\bD_{J_0^c}|_1 \le c_0 |\bD_{J_0}|_1
\end{equation}
we have
\begin{equation*}
\frac{1}{\sqrt{n}}|P_{01} X \bD|_2 \ge  \kappa_1 (s,c_0)
|\bD_{J_{01}}|_2
\end{equation*}
where $P_{01}$ is the projector in $\RR^M$ on the linear span of the
columns of $X_{J_{01}}$.

(ii) Let Assumption 2 be satisfied. Then Assumptions \RE s{c_0} and
\REm sm{c_0} hold with $\kappa (s,c_0) =\kappa (s,m,c_0)=\kappa_2
(s,m,c_0)$. Moreover, for any subset $J_0$ of $\{1,\dots,M\}$ with
cardinality $|J_0|\le s$, and any $\bD\in \RR^M$ such that
(\ref{l21}) holds we have
\begin{equation*}
\frac{1}{\sqrt{n}}|P_{01} X \bD|_2 \ge  \kappa_2 (s,m,c_0)
|\bD_{J_{01}}|_2.
\end{equation*}
\end{lemma}

The proof of the lemma is given in Appendix A.

There exist other sufficient conditions for Assumptions \RE s{c_0}
and \REm sm{c_0} to hold. We mention here three of them implying
Assumption \RE s{c_0}. The first one is the following \cite{b07}.

\noindent {\bf Assumption 3.} {\it For an integer $s$ such that
$1\le s \le M$ we have
$$
\phi_{\min}(s) > 2c_0 \theta_{s,1} \sqrt{s}
$$
where $c_0>0$ is a constant.}

To argue that Assumption 3 implies \RE s{c_0} it suffices to remark
that
\begin{equation*}\begin{split}
\frac{1}{{n}}|X \bD|_2^2 &\ge\frac{1}{n} \bD_{J_0}^T X^T X
\bD_{J_{0}} - \frac{2}{n}|\bD_{J_0}^T X^T
X \bD_{J_{0}^c}|\\
&\ge \phi_{\min}(s)|\bD_{J_{0}}|_2^2 - \frac{2}{n}|\bD_{J_0}^T X^T X
\bD_{J_{0}^c}|
\end{split}\end{equation*}
and, if (\ref{l21}) holds,
\begin{eqnarray*}\label{dla4}
|\bD_{J_0}^T X^T  X \bD_{J_{0}^c}|/n&\le |\bD_{J_{0}^c}|_1
\,\max_{j\in J_0^c}|\bD_{J_{0}}^TX^T {\bf x}_{(j)}|/n
\\
&\le \theta_{s,1}|\bD_{J_{0}^c}|_1|\bD_{J_{0}}|_2
\\
&\le c_0\theta_{s,1}\sqrt{s}|\bD_{J_{0}}|_2^2.
\end{eqnarray*}

Another type of assumption related to ``mutual coherence"
\cite{det04} is discussed in the connection to Lasso in \cite{btw07,
btw06b}. We state it in two different forms given below.

\noindent {\bf Assumption 4.} {\it For an integer $s$ such that
$1\le s \le M$ we have
\begin{equation*}\begin{split}
\phi_{\min}(s) > 2c_0 \theta_{1,1} s
\end{split}\end{equation*}
where $c_0>0$ is a constant.}

It is easy to see that Assumption 4 implies \RE s{c_0}. Indeed, if
(\ref{l21}) holds,
\begin{equation}\begin{split}\label{dla7}
\frac{1}{{n}}|X \bD|_2^2 &\ge\frac{1}{n} \bD_{J_0}^T X^T X
\bD_{J_{0}} - 2\theta_{1,1} |\bD_{J_{0}^c}|_1|\bD_{J_{0}}|_1
\\
&\ge \phi_{\min}(s)|\bD_{J_{0}}|_2^2 - 2c_0 \theta_{1,1}
|\bD_{J_{0}}|_1^2
\\
&\ge (\phi_{\min}(s) - 2c_0 \theta_{1,1} s)|\bD_{J_{0}}|_2^2.
\end{split}\end{equation}
If all the diagonal elements of matrix $X^TX/n$ are equal to 1 (and
thus $\theta_{1,1}$ coincides with the mutual coherence
\cite{det04}), a simple sufficient condition for Assumption \RE
s{c_0} to hold is stated as follows.

\noindent {\bf Assumption 5.} {\it All the diagonal elements of the
Gram matrix $\Psi_n$ are equal to 1 and for an integer $s$ such that
$1\le s \le M$ we have
\begin{equation}\begin{split}\label{dla8}
\theta_{1,1} < \frac{1}{ (1+2c_0)  s}\,.
\end{split}\end{equation}
where $c_0>0$ is a constant.}

In fact, separating the diagonal and off-diagonal terms of the
quadratic form we get
\begin{eqnarray*}
\bD_{J_0}^T X^T  X \bD_{J_{0}}/n \ge |\bD_{J_{0}}|_2^2 -
\theta_{1,1} |\bD_{J_{0}}|_1^2\ge |\bD_{J_{0}}|_2^2 (1 -
\theta_{1,1} s).
\end{eqnarray*}
Combining this inequality with (\ref{dla7}) we see that Assumption
\RE s{c_0} is satisfied whenever (\ref{dla8}) holds.

Unfortunately, Assumption \RE s{c_0} has some weakness. Let, for
example, $f_j$, $j=1,\dots,2^m-1$, be the Haar wavelet basis on
$[0,1]$ ($M=2^m$) and consider $Z_i=i/n$, $i=1,\dots,n$. If $M\gg
n$, it is clear that $\phi_{\min}(1)=0$ since there are functions
$f_j$ on the highest resolution level whose supports (of length
$M^{-1}$) contain no points $Z_i$. So, none of the Assumptions 1 --
4 holds. A less severe although similar situation is when we
consider step functions: $f_j(t)=I_{\{t<j/M\}}$ for $t\in[0,1]$. It
is clear that $\phi_{\min}(2)=O(1/M)$, although sparse
representation in this basis is very natural. Intuitively, the
problem arises only because we include very high resolution
components. Therefore, we may try to restrict the set $J_0$ in \RE
s{c_0} to low resolution components, which is quite reasonable
because the ``true" or ``interesting" vectors of parameters $\beta$
are often characterized by such $J_0$. This idea is formalized in
Section \ref{sec5}, cf. Corollary \ref{cor2}, see also a remark
after Theorem \ref{th5} in Section \ref{sec6}.

\section{Approximate equivalence}
\label{sec3}

In this section we prove a type of approximate equivalence between the Lasso
and the Dantzig selector. It is expressed as closeness of the prediction
losses $\| \widehat f_D -f\|_n^2$ and $\| \wideth{f}_L -f \|_n^2$ when the
number of non-zero components of the Lasso or the Dantzig selector  is small
as compared to the sample size.
\begin{thm}\label{th1}
Let  $W_i$ be independent ${\mathcal N}(0,\sigma^2)$ random
variables with $\sigma^2>0$. Fix $n\ge1$, $M\ge2$. Let Assumption
\RE s1 be satisfied with $1\le s\le M$. Consider the Dantzig
estimator $\wideth{f}_D$ defined by (\ref{da1}) -- (\ref{da2}) with
$$
r= A\sigma \sqrt{\frac{\log M}{n}}\,
$$
where $A>2\sqrt 2$, and the Lasso estimator $\wideth{f}_L$ defined
by (\ref{lass}) -- (\ref{lse}) with the same $r$.

If $\scm(\wh\beta_L)\le s$, then with probability at least
$1-M^{1-A^2/8}$ we have
\begin{equation}\begin{split}\label{th11}
&\Bigl| \| \wideth f_D -f\|_n^2 - \| \wideth{f}_L -f \|_n^2\Bigr| \leq
16 A^2 \frac{\scm(\wh\beta_L)\sigma^2}{n}\frac{f_{\max}^2}{\kappa^2(s,1)}
\log M.
\end{split}\end{equation}

\end{thm}

Note that the RHS of \eqref{th11}  is bounded by a product of three
factors (and a numerical constant which, unfortunately, equals at
least 128). The first factor, $\scm(\wh\beta_L)\sigma^2/n \le
s\sigma^2/n$, corresponds to the error rate for prediction in
regression with $s$ parameters. The two other factors, $\log M$ and
$f^2_{\max}/\kappa^2(s,1)$, can be regarded as a price to pay for
the large number of regressors. If the Gram matrix $\Psi_n$ equals
the identity matrix (the white noise model), then there is only the
$\log M$ factor. In the general case, there is another factor,
$f^2_{\max}/\kappa^2(s,1)$ representing the extent to which the Gram
matrix is ill-posed for estimation of sparse vectors.

We also have the following result that we state for simplicity under
the assumption that $\|f_j\|_n=1,\ j=1,\ldots,M$. It gives a bound
in the spirit of Theorem \ref{th1} but with $\scm(\wh\beta_D)$
rather than $\scm(\wh\beta_L)$ on the right hand side.

\begin{thm}\label{th2}
Let the assumptions of Theorem \ref{th1} hold, but with \RE s5 in
place of \RE s1, and let $\|f_j\|_n=1,\ j=1,\ldots,M$. If
$\scm(\wh\beta_D)\le s$, then with probability at least
$1-M^{1-A^2/8}$ we have
\begin{equation}\begin{split} &\| \wideth f_L -f\|_n^2 \le 10\| \wideth{f}_D -f \|_n^2 +
 81A^2  \frac{\scm(\wh\beta_D)
 \sigma^2}{n}\frac{\log M}{\kappa^2(s,5)}.
\end{split}\end{equation}
\end{thm}

\noindent{\sc Remark.} The approximate equivalence is essentially
that of the rates as Theorem \ref{th1} exhibits. A statement free of
$\scm(\beta)$ holds for linear regression, see discussion after
Theorem \ref{th5} and Theorem \ref{th4a} below.

\section{Oracle inequalities for prediction loss}
\label{sec5}

Here we prove sparsity oracle inequalities for the prediction loss
of Lasso and Dantzig estimators. These inequalities allow us to
bound the difference between the prediction errors of the estimators
and the best sparse approximation of the regression function (by an
oracle that knows the truth, but is constrained by sparsity). The
results of this section, together with those of Section \ref{sec3},
show that the distance between the prediction losses of Dantzig and
Lasso estimators is of the same order as the distances between them
and their oracle approximations.

A general discussion of sparsity oracle inequalities can be found in \cite{t06}. Such inequalities
have been recently obtained for the Lasso type estimators in a
number of settings
\cite{btw04,btw06a,btw07,btw06b,btw07a,k06,vdg06}. In particular,
the regression model with fixed design that we study here is
considered in \cite{btw04,btw06a,btw07}. The assumptions on the Gram
matrix $\Psi_n$ in \cite{btw04,btw06a,btw07} are more restrictive
than ours: in those papers either $\Psi_n$ is positive definite or a
mutual coherence condition similar to (\ref{dla8}) is imposed.

\begin{thm}\label{th-ora-las}
Let  $W_i$ be independent ${\mathcal N}(0,\sigma^2)$ random
variables with $\sigma^2>0$. Fix some $\eps>0$ and integers $n\ge1$,
$M\ge2$, $1\le s\le M$. Let Assumption \RE
s{(3+4/\eps)f_{\max}/f_{\min}} be satisfied. Consider the Lasso
estimator $\wideth{f}_L$ defined by (\ref{lass}) -- (\ref{lse}) with
$$
r= A\sigma \sqrt{\frac{\log M}{n}}\,
$$
for some $A>2\sqrt{2}$. Then, with probability at least
$1-M^{1-A^2/8}$, we have
\begin{equation}\label{predlasso}
\begin{split}
&\hspace{-0.3em} \| \wideth f_L -f\|_n^2
\\
&\le  (1+\eps) \inf_{\substack{\beta\in\RR^M:\\\scm(\beta)\le
s}}\left\{\| {\sf f}_\beta -f \|_n^2 + \,\frac{C(\eps)f_{\max}^2A^2
\sigma^2}{\kappa^2}\, \frac{\scm(\beta)\log M}{n}\right\}\,
\end{split}\end{equation}
where $\kappa=\kappa(s,(3+4/\eps)f_{\max}/f_{\min})$, and
$C(\eps)>0$ is a constant depending only on $\eps$.
\end{thm}

We now state as a corollary a softer version of Theorem
\ref{th-ora-las} that can be used to eliminate the pathologies
mentioned at the end of Section \ref{discRE}. For this purpose we
define
$${\mathcal J}_{s,\gamma,c_0} = \Big\{J_0 \subset\{1,\dots,M\}: |J_0|\leq
s \ \ \text{and} \ 
\min_{\substack{\bD\neq 0, \\ |\bD_{J_0^c}|_1\le c_0|\bD_{J_0}|_1}}
\frac{|X\bD|_2}{\sqrt{n}|\bD_{J_0}|_2}\ge \gamma\Big\}$$ where
$\gamma>0$ is a constant, and set
$$
\Lambda_{s,\gamma,c_0}= \{\beta:
J(\beta)\in {\mathcal J} _{s,\gamma,c_0}\}.
$$
In similar way, we define ${\mathcal J}_{s,\gamma,m, c_0}$ and
$\Lambda_{s,\gamma,m,c_0}$ corresponding to Assumption
\REm sm{c_0}.

\begin{cor}\label{cor2} Let  $W_i$, $s$ and the Lasso estimator
$\wideth{f}_L$ be the same as in Theorem \ref{th-ora-las}. Then,
for all $n\ge1$,  $\eps>0$, and $\gamma>0$, with probability at least
$1-M^{1-A^2/8}$ we have
\begin{equation*}
\begin{split} &\hspace{-1em}\| \wideth f_L -f\|_n^2
\\
& \le (1+\eps)
\inf_{\beta\in\bar\Lambda_{s,\gamma,\eps}}\left\{\| {\sf
f}_\beta -f \|_n^2 + \,\frac{C(\eps)f_{\max}^2A^2
\sigma^2}{\gamma^2}\, \left(\frac{\scm(\beta)\log
M}{n}\,\right)\right\}\,
\end{split}
\end{equation*}
where $\bar\Lambda_{s,\gamma,\eps}
=\{\beta\in\Lambda_{s,\gamma,(3+4/\eps)f_{\max}/f_{\min}}:\,\scm(\beta)\le
s\}$.
\end{cor}
To obtain this corollary it suffices to observe that the proof of
Theorem \ref{th-ora-las} goes through if we drop Assumption \RE
s{(3+4/\eps)f_{\max}/f_{\min}} but we assume instead that
$\beta\in\Lambda_{s,\gamma,(3+4/\eps)f_{\max}/f_{\min} }$ and we
replace $\kappa(s,(3+4/\eps)f_{\max}/f_{\min})$ by $\gamma$.

We would like now to get a sparsity oracle inequality similar to
that of Theorem \ref{th-ora-las} for the Dantzig estimator
$\wideth{f}_D$. We will need a mild additional assumption on $f$.
This is due to the fact that not every $\beta\in\RR^M$ obeys to
the Dantzig constraint, and thus we cannot assure the key relation
(\ref{ll22}) for all $\beta\in\RR^M$. One possibility would be to
prove inequality as (\ref{predlasso}) where the infimum on the right
hand side is taken over $\beta$ satisfying not only $\scm(\beta)\le
s$ but also the Dantzig constraint. However, this seems not very
intuitive since we cannot guarantee that the corresponding ${\sf
f}_\beta$ gives a good approximation of the unknown function $f$.
Therefore we choose another approach (cf. \cite{btw06b}): we
consider $f$ satisfying the {\it weak sparsity} property relative to
the dictionary $f_1,\dots,f_M$. That is, we assume that there exist
an integer $s$ and constant $C_0<\infty$ such that the set
\begin{equation}\label{La}
\Lambda_s = \left\{ \beta \in \RR^M: \ \scm(\beta)\le s,  \ \| {\sf
f}_\beta - f\|_n^2 \le \frac{C_0 f_{\max}^2r^2}{\kappa^2}
\scm(\beta) \right\}
\end{equation}
is non-empty. Here $\kappa$ is the same as in Theorem
\ref{th-ora-las}. The second inequality in (\ref{La}) says that the
``bias" term $\| {\sf f}_\beta - f\|_n^2$ cannot be much larger than
the ``variance term" $\sim{f_{\max}^2r^2}\kappa^{-2} \scm(\beta)$,
cf. (\ref{predlasso}). Weak sparsity is milder than the sparsity
property in the usual sense: the latter means that $f$ admits the
exact representation $f={\sf f}_{\beta^*}$ for some $\beta^*\in
\RR^M$, with hopefully small $\scm(\beta^*)= s$.\\

\begin{proposition}\label{th-ora-dan}
Let  $W_i$ be independent ${\mathcal N}(0,\sigma^2)$ random
variables with $\sigma^2>0$. Fix some $\eps>0$ and integers $n\ge1$,
$M\ge2$. Let $f$ obey the weak sparsity assumption for some
$C_0<\infty$ and some $s$ such that $1\le  s\max \{C_1(\eps),1\} \le
M$ where
$$
C_1(\eps)=4\left[(1+\eps)C_0 + C(\eps)\right]
\frac{\phi_{\max}f_{\max}^2}{\kappa^2 f_{\min}^{2}}\,
$$
and $C(\eps)$, $\kappa$ are the constants in Theorem
\ref{th-ora-las}. Suppose further that Assumption \RE {s\max
\{C_1(\eps),1\}}{(3+4/\eps)f_{\max}/f_{\min}} is satisfied. Consider
the Dantzig estimator $\wideth{f}_D$ defined by (\ref{da1}) --
(\ref{da2}) with
$$
r= A\sigma \sqrt{\frac{\log M}{n}}\,
$$
and $A>2\sqrt{2}$. Then, with probability at least $1-M^{1-A^2/8}$,
we have
\begin{equation}\label{preddan}
\begin{split}
&\hspace{-1em}\| \wideth f_D -f\|_n^2
\\
&\le (1+\eps)
\inf_{\beta\in\RR^M:\,\scm(\beta)=s}\| {\sf f}_\beta -f \|_n^2 +
C_2(\eps)\,\frac{f_{\max}^2A^2 \sigma^2}{\kappa_0^2}\,
\left(\frac{s\log M}{n}\,\right)\,.
\end{split}
\end{equation}
Here $C_2(\eps)=16 C_1(\eps)+ C(\eps)$ and $\kappa_0= \kappa(\max
(C_1(\eps),1)s, (3+4/\eps)f_{\max}/f_{\min})$.
\end{proposition}

Note that the sparsity oracle inequality (\ref{preddan}) is slightly
weaker than the analogous inequality (\ref{predlasso}) for the
Lasso: we have here $\inf_{\beta\in\RR^M:\,\scm(\beta)=s}$ instead
of $\inf_{\beta\in\RR^M:\,\scm(\beta)\le s}$ in (\ref{predlasso}).

\section{Special case: parametric estimation in linear regression}
\label{sec6}

In this section we assume that the vector of observations
$\bmy=(Y_1,\dots,Y_n)^T$ is of the form
\begin{equation}\begin{split}\label{linmod}
\bmy=X\bb^*+ \bmw
\end{split}\end{equation}
where $X$ is an $n\times M$ deterministic matrix, $\bb^*\in\RR^M$ and
$\bmw=(W_1,\dots,W_n)^T$.

We consider dimension $M$ that can be of order $n$ and even much
larger. Then $\bb^*$ is, in general, not uniquely defined. For
$M>n$, if (\ref{linmod}) is satisfied for $\bb^*=\bb_0$ there exists
an affine space ${\cal U}=\{\bb^*:X\bb^* = X\bb_0\}$ of vectors
satisfying (\ref{linmod}). The results of this section are valid for
any $\bb^*$ such that (\ref{linmod}) holds. However, we will suppose
that Assumption \RE s{c_0} holds with $c_0\ge1$ and that
$\scm(\beta^*)\le s$. Then the set ${\cal U}\cap \{\bb^*:
\scm(\beta^*)\le s\}$ reduces to a single element (cf. Remark 2 at
the end of this section). In this sense, there is a unique sparse
solution of (\ref{linmod}).

Our goal in this section, unlike that of the previous ones, is to
estimate both $X\bb^*$ for the purpose of prediction and $\bb^*$
itself for purpose of model selection. We will see that meaningful
results are obtained when the sparsity index $\scm(\bb^*)$ is small.

It will be assumed throughout this section that the diagonal
elements of the Gram matrix $\Psi_n = X^TX/n$ are all equal to 1
(this is equivalent to the condition $\|f_j\|_n=1,\ j=1,\ldots,M,$
in the notation of previous sections). Then the Lasso estimator of
$\bb^*$ in (\ref{linmod}) is defined by
\begin{equation}\begin{split}\label{Lasso}
\wh{\bb}_L = \argmin_{\bb\in\RR^M} \left\{ \frac{1}{n}|\bmy-X\bb|_2^2 +
2r|\bb|_1 \right\}.
\end{split}\end{equation}
The correspondence between the notation here and that of the previous sections is the following:
$$
\| {\sf f}_{\beta}\|_n^2 = |X\,\bb|_2^2/n, \quad \| {\sf
f}_{\beta} - f\|_n^2 = |X\,(\bb-\bb^*)|_2^2/n, \quad \| \wideth f_L - f\|_n^2=|X\,(\widehat\bb_L-\bb^*)|_2^2/n.
$$
The Dantzig selector for linear model (\ref{linmod}) is defined by
\begin{equation}\begin{split}\label{dl4}
\widehat \bb_D = \argmin_{\bb\in\Lambda}|\bb|_1
\end{split}\end{equation}
where
$$
\Lambda= \Big\{\bb\in\RR^M:\ \Big|\frac{1}{n}X^T(\bmy-X\bb)\Big|_\infty
\le r \Big\}
$$
is the set of all $\bb$ satisfying the Dantzig constraint.

We first get bounds on the rate of convergence of Dantzig selector.

\begin{thm}\label{th4}
Let  $W_i$ be independent ${\mathcal N}(0,\sigma^2)$ random
variables with $\sigma^2>0$, let all the diagonal elements of the
matrix $X^TX/n$ be equal to 1, and $\scm(\bb^*)\le s$, where $1\le s
\le M$, $n\ge1$, $M\ge2$. Let Assumption \RE s1 be satisfied.
Consider the Dantzig selector $\widehat \bb_D$ defined by
(\ref{dl4}) with
$$
r= A\sigma \sqrt{\frac{\log M}{n}}\,
$$
and $A>\sqrt{2}$. Then, with probability at least $1-M^{1-A^2/2}$,
we have
\begin{align}\label{th43}
 |\widehat\bb_D-\bb^*|_1 &\le \frac{8A}{\kappa^2(s,1)}
\sigma \, s \sqrt{\frac{\log M}{n}},\\
\label{th44} |X(\widehat\bb_D-\bb^*)|_2^2 &\le
\frac{16A^2}{\kappa^2(s,1)}\, \sigma^2\, s\log M.
\end{align}
If Assumption \REm sm1 is satisfied, then with the same probability
as above, simultaneously for all $1<p\le2$ we have
\begin{equation}\begin{split}
\label{th42} &|\widehat\bb_D-\bb^*|_p^p \le 2^{p-1}8
\left\{1+\sqrt{\frac{s}{m}}\,\right\}^{2(p-1)} \, s \, \left(
\frac{A\sigma}{\kappa^2(s,m,1)}\, \sqrt{\frac{ \log M}{n}}\right)^p.
\end{split}\end{equation}

\end{thm}

Note that, since $s\le m$, the factor in curly brackets in
(\ref{th42}) is bounded by a constant independent of $s$ and $m$.
Under Assumption  1 in Section \ref{discRE} with $c_0=1$ (which is
less general than \REm ss1, cf.
 Lemma \ref{l2}(i)) a bound
of the form (\ref{th42}) for the case $p=2$ is established by Candes
and Tao \cite{ct07}.

Bounds on the rate of convergence of the Lasso selector are quite
similar to those obtained in Theorem \ref{th4}. They are given by
the following result.

\begin{thm}\label{th5}
Let  $W_i$ be independent ${\mathcal N}(0,\sigma^2)$ random
variables with $\sigma^2>0$. Let all the diagonal elements of the
matrix $X^TX/n$ be equal to 1, and $\scm(\bb^*)\le s$ where $1\le s
\le M$, $n\ge1$, $M\ge2$. Let Assumption \RE s3 be satisfied.
Consider the Lasso estimator $\widehat \bb_L$ defined by
(\ref{Lasso}) with
$$
r= A\sigma \sqrt{\frac{\log M}{n}}\,
$$
and $A>2\sqrt{2}$. Then, with probability at least $1-M^{1-A^2/8}$,
we have
\begin{align}
\label{th53} &|\widehat\bb_L-\bb^*|_1 \le \frac{16A}{\kappa^2(s,3)}
\sigma \, s \sqrt{\frac{\log M}{n}},\\
\label{th54} &|X(\widehat\bb_L-\bb^*)|_2^2 \le
\frac{16A^2}{\kappa^2(s,3)}\, \sigma^2\, s\log M,
\\
\label{th55} &\scm(\widehat\bb_L) \le \frac{64 \phi_{\max}}{\kappa^2(s,3)}\,s.
\end{align}
If Assumption \REm sm3 is satisfied, then with the same probability
as above, simultaneously for all $1<p\le2$ we have
\begin{equation}\begin{split}
\label{th52} &|\widehat\bb_L-\bb^*|_p^p \le 16
\left\{1+3\sqrt{\frac{s}{m}}\,\right\}^{2(p-1)}\, s\,\left(
\frac{A\sigma }{\kappa^2(s,m,3)}\, \sqrt{\frac{ \log M}{n}}\right)^p.
\end{split}\end{equation}

\end{thm}

Inequalities of the form similar to (\ref{th53}) and (\ref{th54})
can be deduced from the results of \cite{btw06a} under more
restrictive conditions on the Gram matrix (the mutual coherence
assumption, cf. Assumption 5 of Section \ref{discRE}).

Assumptions \RE s1 respectively \RE s3 can be dropped in
Theorem \ref{th4} and \ref{th5} if we assume $\beta^*\in
\Lambda_{s,\gamma,c_0}$ with $c_0=1$ or $c_0=3$ as appropriate. Then
\eqref{th43}, \eqref{th44} or respectively \eqref{th53},
\eqref{th54} hold with $\kappa=\gamma$. This is analogous to
Corollary \ref{cor2}. Similarly \eqref{th42} and \eqref{th52} hold
with $\kappa=\gamma$ if $\beta^*\in \Lambda_{s, \gamma,m,c_0}$ with
$c_0=1$ or $c_0=3$ as appropriate.

Observe that combining Theorems \ref{th4} and \ref{th5} we can
immediately get bounds for the differences between Lasso and Dantzig
selector $|\widehat\bb_L-\widehat\bb_D|_p^p$ and
$|X(\widehat\bb_L-\widehat\bb_D)|_2^2$. Such bounds have the same
form as those of Theorems \ref{th4} and \ref{th5}, up to numerical
constants. Another way of estimating these differences follows
directly from the proof of Theorem \ref{th4}. It suffices to observe
that the only property of $\bb^*$ used in that proof is the fact
that $\bb^*$ satisfies the Dantzig constraint on the event of given
probability, which is also true for the Lasso solution $\wh\bb_L$.
So, we can replace $\bb^*$ by $\wh\bb_L$ and $s$ by $\scm(\wh\bb_L)$
everywhere in Theorem \ref{th4}. Generalizing a bit more, we easily
derive the following fact.
\begin{thm}\label{th4a}
The result of Theorem \ref{th4} remains valid if we replace there
$|\widehat\bb_D-\bb^*|_p^p$ by $\sup\{|\widehat\bb_D-\bb|_p^p:\,
\bb\in\Lambda, \scm(\bb)\le s\}$ for $1\le p\le 2$ and
$|X(\widehat\bb_D-\bb^*)|_2^2$ by
$\sup\{|X(\widehat\bb_D-\bb)|_2^2:\, \bb\in\Lambda, \scm(\bb)\le
s\}$ respectively. Here $\Lambda$ is the set of all vectors
satisfying the Dantzig constraint.
\end{thm}

\noindent{\sc Remarks.}

\begin{list}{}
{\renewcommand{\makelabel}{\addtocounter{enumi}{1}\arabic{enumi}.}
\setlength{\labelwidth}{1em}
\setlength{\itemindent}{0.51em}
\setlength{\itemsep}{2ex}
\setlength{\leftmargin}{0.15cm}
\setlength{\listparindent}{\parindent}}

\item
Theorems \ref{th4} and \ref{th5} only give non-asymptotic upper
bounds on the loss, with some probability and under some conditions.
The probability depends on $M$ and the conditions depend on $n$ and
$M$: recall that Assumptions RE($s,c_0$) and RE($s,m,c_0$) are
imposed on the $n\times M$ matrix $X$. To deduce asymptotic
convergence (as $n\to \infty$ and/or as $M\to\infty$) from Theorems
\ref{th4} and \ref{th5} we would need some very strong additional
properties, such as simultaneous validity
 of Assumption RE($s,c_0$) or
RE($s,m,c_0$) (with one and the same constant $\kappa$) for
infinitely many $n$ and $M$.

\item

Note that Assumptions RE($s,c_0$) or RE($s,m,c_0$) do not imply
identifiability of $\beta^*$ in the linear model (\ref{linmod}).
However, the vector $\beta^*$ appearing in the statements of
Theorems \ref{th4} and \ref{th5} is uniquely defined because we
suppose there in addition that $\scm(\beta^*)\le s$ and $c_0\ge1$.
Indeed, if there exists a $\beta'$ such that $X\beta'=X\beta^*$, and
$\scm(\beta')\le s$ then in view of assumption \RE {s}{c_0} with
$c_0\geq 1$ we have necessarily $\beta^*=\beta'$ (cf. discussion
following the definition of \RE s{c_0}). On the other hand, Theorem
\ref{th4a} applies to certain values of $\bb$ that do not come from
the model (\ref{linmod}) at all.

\item
For the smallest value of $A$ (which is $A=2\sqrt{2}$) the constants
in the bound of Theorem \ref{th5} for the Lasso are larger than the
corresponding numerical constants for the Dantzig selector given in
Theorem \ref{th4}, again for the smallest admissible value
$A=\sqrt{2}$. On the contrary, the Dantzig selector has certain
defects as compared to Lasso when the model is nonparametric, as
discussed in Section \ref{sec5}. In particular, to obtain sparsity
oracle inequalities for the Dantzig selector we need some
restrictions on $f$, for example the weak sparsity property. On the
other hand, the sparsity oracle inequality (\ref{predlasso}) for the
Lasso is valid with no restriction on $f$.

\item
The proofs of Theorems \ref{th4} and \ref{th5} differ mainly in the
value of the tuning constant: $c_0=1$ in Theorem \ref{th4} and
$c_0=3$ in Theorem \ref{th5}. Note that since the Lasso solution
satisfies the Dantzig constraint we could have obtained a result
similar to Theorem \ref{th5}, though with less accurate numerical
constants, by simply conducting the proof of Theorem \ref{th4} with
$c_0=3$. However, we act differently: we deduce (\ref{dl8a})
directly from (\ref{onA}), and not from (\ref{dl5}). This is done
only for the sake of improving the constants: in fact, using
(\ref{dl5}) with $c_0=3$ would yield (\ref{dl8a}) with the doubled
constant on the right hand side.

\item
For the Dantzig selector in the linear regression model and under
Assumptions 1 or 2 some further improvement of constants in the
$\ell_p$ bounds for the coefficients can be achieved by applying the
general version of Lemma \ref{l2} with the projector $P_{01}$
inside. We do not pursue this issue here.

\item
All our results are stated with probabilities at least
$1-M^{1-A^2/2}$ or $1-M^{1-A^2/8}$. These are reasonable (but not
the most accurate) lower bounds on the probabilities $\PP(\B)$ and
$\PP(\A)$ respectively: we have chosen them just for readability.
 Inspection of (\ref{tail}) shows that
they can be refined to $1-2M\Phi(A\sqrt{\log M})$ and
$1-2M\Phi(A\sqrt{\log M}/2)$ respectively where $\Phi(\cdot)$ is the
standard normal c.d.f.
\end{list}
\appendix
\section{}
\label{app:re}

{\sc Proof of Lemma \ref{l2}}. Consider a partition $J_0^c$ into
subsets of size $m$, with the last subset of size $\le m$:
$J_0^c=\cup_{k=1}^K J_k$ where $K\ge 1$, $|J_k|=m$ for
$k=1,\dots,K-1$ and $|J_K|\le m$, such that $J_k$ is the set of
indices corresponding to $m$ largest in absolute value coordinates
of $\bD$ outside $\cup_{j=1}^{k-1} J_j$ (for $k<K$) and $J_K$ is the
remaining subset. We have
\begin{equation}\begin{split}\label{dl1}
|P_{01} X \bD|_2&\ge |P_{01} X \bD_{J_{01}}|_2 - \Big|\sum_{k=2}^K
P_{01} X \bD_{J_{k}}\Big|_2
\\
&= |X \bD_{J_{01}}|_2 - \Big|\sum_{k=2}^K P_{01} X
\bD_{J_{k}}\Big|_2\\  &\ge|X \bD_{J_{01}}|_2 - \sum_{k=2}^K |P_{01}
X \bD_{J_{k}}|_2.
\end{split}\end{equation}
We will prove first part (ii) of the lemma. Since for $k\ge1$ the
vector $\bD_{J_{k}}$ has only $m$ non-zero components we obtain
\begin{equation}\begin{split}\label{dl*}
\frac{1}{\sqrt{n}}|P_{01} X \bD_{J_{k}}|_2 \le\frac{1}{\sqrt{n}}| X
\bD_{J_{k}}|_2 \le \sqrt{\phi_{\max}(m)}\,|\bD_{J_{k}}|_2.
\end{split}\end{equation}
 Next, as in \cite{ct07}, we observe
that $|\bD_{J_{k+1}}|_2\le |\bD_{J_{k}}|_1/\sqrt{m}$,
$k=1,\dots,K-1$, and therefore
\begin{equation}\label{dl2}
\begin{split}
\sum_{k=2}^K |\bD_{J_{k}}|_2 \le \frac{|\bD_{J_{0}^c}|_1}{\sqrt{m}}
\le \frac{c_0|\bD_{J_{0}}|_1}{\sqrt{m}} \le c_0 \sqrt{\frac{s}{m}}\
|\bD_{J_{0}}|_2 \le c_0 \sqrt{\frac{s}{m}}\ |\bD_{J_{01}}|_2
\end{split}
\end{equation}
where we used (\ref{l21}). From (\ref{dl1}) -- (\ref{dl2}) we find
\begin{equation*}
\begin{split}
\frac{1}{\sqrt{n}}|X \bD|_2&\ge  \frac{1}{\sqrt{n}}|X
\bD_{J_{01}}|_2 - c_0 \sqrt{\phi_{\max}(m)}\sqrt{\frac{s}{m}}\
|\bD_{J_{01}}|_2
\\
&\ge \left(\sqrt{\phi_{\min}(s+m)}-c_0
\sqrt{\phi_{\max}(m)}\sqrt{\frac{s}{m}}\,\right)|\bD_{J_{01}}|_2
\end{split}
\end{equation*}
which proves part (ii) of the lemma.

The proof of part (i) is analogous. The only difference is that we
replace in the above argument $m$ by $s$ and instead of (\ref{dl*})
we use the following bound (cf. \cite{ct07}):
\begin{equation*}\begin{split}
\frac{1}{\sqrt{n}}|P_{01} X \bD_{J_{k}}|_2 \le
\frac{\theta_{s,2s}}{\sqrt{\phi_{\min}(2s)}}\,|\bD_{J_{k}}|_2.
\end{split}\end{equation*}

\section{Two lemmata and the proofs of the results}
\label{app:proofs}
\begin{lemma}\label{ll1} Fix $M\ge2$ and $n\ge1$.
Let  $W_i$ be independent ${\mathcal N}(0,\sigma^2)$ random
variables with $\sigma^2>0$ and let $\wideth f_L$ be the Lasso
estimator defined by (\ref{lse}) with
$$
r= A\sigma \sqrt{\frac{\log M}{n}},
$$
for some $A>2\sqrt{2}$. Then, with probability at least
$1-M^{1-A^2/8}$ we have simultaneously for all $\beta\in\RR^M$:
\begin{equation}\begin{split}\label{onA}
&\hspace{-3em}\| \widehat f_L - f\|_n^2 + r\sum_{j=1}^M \|f_j\|_{n}     |\widehat \beta_{j,L}
- \beta_j|
\\
&\le \| {\sf f}_{\beta} - f\|_n^2 + 4r \sum_{j\in
J(\beta)} \|f_j\|_{n}|\widehat \beta_{j,L} - \beta_j|
\\
&\le \| {\sf f}_{\beta} - f\|_n^2 +
4r \sqrt{\scm(\beta)}\sqrt{\sum_{j\in J(\beta)} \|f_j\|_{n}^2|\widehat
\beta_{j,L} - \beta_j|^2},
\end{split}\end{equation}
and
\begin{equation}\begin{split}\label{dcf}
\Big|\frac{1}{n}X^T({\bmf }-X\widehat\beta_L)\Big|_\infty \le
3rf_{\max}/2.
\end{split}\end{equation}
Furthermore, with the same probability
\begin{equation}\begin{split}\label{Lassom}
\scm(\wh\beta_L) \le 4 \phi_{\max} f_{\min}^{-2} \left(\| \widehat f_L
- f\|_n^2/r^2 \right)
\end{split}\end{equation}
where $\phi_{\max}$ denotes the maximal eigenvalue of the matrix
$X^TX/n$.
\end{lemma}

\begin{proof}[Proof of Lemma \ref{ll1}] The result (\ref{onA}) is essentially Lemma 1
from \cite{btw06b}. For completeness, we give its proof. Set $
r_{n,j} = r \|f_j\|_n $. By definition,
\begin{eqnarray*}
\wh{S}(\widehat\beta_L) +  2\sum_{j=1}^M r_{n,j} |\widehat
\beta_{j,L}| \le \wh{S}(\beta) + 2\sum_{j=1}^M r_{n,j} |\beta_j|
\end{eqnarray*} for all $\beta\in\RR^M$,
which is equivalent to
\begin{eqnarray*}
\| \wideth{f}_L - f\|_n^2 +  2\sum_{j=1}^M r_{n,j} |\widehat
\beta_{j,L}| \le \| {\sf f}_{\beta} - f\|_n^2 +  2\sum_{j=1}^M
r_{n,j} |\beta_j| +\frac{2}{n}\sum_{i=1}^n W_i (\wideth{f}_L -
{\sf f}_{\beta})(Z_i) .
\end{eqnarray*}
Define the random variables $ V_j =  n^{-1} \sum_{i=1}^n f_j(Z_i)
W_i,\quad 1\le j\le M,$ and the event
\[ \A = \bigcap_{j=1}^M \left\{  2| V_j | \leq r_{n,j} \right\}.\]
Using an elementary bound on the tails of Gaussian disribution we
find that the probability of the complementary event $\A^c$
satisfies
\begin{equation}\label{tail}
\begin{split}
\PP \{ \A^c\}  &\le \sum_{j=1}^M  \PP\{ \sqrt{n}|V_j| >
\sqrt{n}r_{n,j}/2\} \le M \,\PP\{ |\eta|\ge r \sqrt{n}/(2\sigma)\}
\\
&\le M \exp\left( -\frac{n r^2}{ 8 \sigma^2} \right) = M \exp\left(
-\frac{A^2\log M}{ 8} \right) = M^{1-A^2/8}
\end{split}
 \end{equation}
where $\eta\sim{\mathcal N}(0,1)$. On the event $\A$ we have
\begin{eqnarray*}
\| \wideth f_L - f\|_n^2 & \le&  \| {\sf f}_{\beta} - f\|_n^2 +
\sum_{j=1}^M r_{n,j}|\widehat \beta_{j,L} - \beta_j| + \sum_{j=1}^M
2r_{n,j} |\beta_j| -  \sum_{j=1}^M 2r_{n,j} |\widehat\beta_{j,L}|.
\end{eqnarray*}
Adding the term $\sum_{j=1}^M r_{n,j}|\widehat \beta_{j,L} -
\beta_j|$ to both sides of this inequality yields, on $\A$,
\begin{eqnarray*}
\| \wideth f_L - f\|_n^2 + \sum_{j=1}^M r_{n,j}|\widehat \beta_{j,L}
- \beta_j|  &\le \| {\sf f}_{\beta} - f\|_n^2 + 2\sum_{j=1}^M
r_{n,j}\left(|\widehat \beta_{j,L} - \beta_j| + |\beta_j| -
|\widehat\beta_{j,L}|\right).
\end{eqnarray*}
Now, $|\widehat \beta_{j,L} - \beta_j| + |\beta_j| -
|\widehat\beta_{j,L}|=0$ for $j\not\in  J(\beta)$, so that on $\A$ we
get (\ref{onA}).

To prove (\ref{dcf}) it suffices to note that on $\A$ we have
\begin{equation}\begin{split}\label{dc2}
\Big|\frac{1}{n}D^{-1/2}X^TW\Big|_\infty \le r/2.
\end{split}\end{equation}
Now, $\bmy={\bmf }+ \bmw$, and (\ref{dcf}) follows from (\ref{dc}),
(\ref{dc2}).

 We finally prove (\ref{Lassom}). The necessary and
sufficient condition for $\wh \beta_L$ to be the Lasso solution can
be written in the form
\begin{equation}\label{las1}
\begin{split}
\frac{1}{n}{\mathbf x}_{(j)}^T(y-X\wh \beta_L) &= r \|f_j\|_n
\,{\rm
sign}(\wh \beta_{j,L}) \quad\mbox{if \ $\wh\beta_{j,L}\neq 0$},
\\
\Big|\frac{1}{n}{\mathbf x}_{(j)}^T(y-X\wh \beta_L)\Big| &\le r
\|f_j\|_n \qquad \mbox{if \ $\wh\beta_{j,L}= 0$}
\end{split}
\end{equation}
where ${\mathbf x}_{(j)}$ denotes the $j$th column of $X$,
$j=1,\dots,M$. Next, (\ref{dc2}) yields that on $\A$ we have
\begin{equation}\begin{split}\label{las2}
\Big|\frac{1}{n}{\mathbf x}_{(j)}^TW\Big| &\le r \|f_j\|_n/2, \quad
j=1,\dots,M.
\end{split}\end{equation}
Combining (\ref{las1}) and (\ref{las2}) we get
\begin{equation}\begin{split}\label{las3}
\Big|\frac{1}{n}{\mathbf x}_{(j)}^T({\bmf}-X\wh \beta_L)\Big| &\ge
r \|f_j\|_n/2  \quad\mbox{if \ $\wh\beta_{j,L}\neq 0$}.
\end{split}\end{equation}
Therefore,
\begin{eqnarray*}
\frac{1}{n^2}({\bmf}-X\wh \beta_L)^T XX^T ({\bmf}-X\wh \beta_L)
&= \frac{1}{n^2}\sum_{j=1}^M \Big({\mathbf x}_{(j)}^T({\bmf}-X\wh
\beta_L)
\Big)^2\\
&\ge \frac{1}{n^2}\sum_{j:\,\wh\beta_{j,L}\neq 0} \Big({\mathbf
x}_{(j)}^T({\bmf}-X\wh \beta_L)\Big)^2 \\
&= \scm(\wh\beta_L) r^2 \|f_j\|_n^2/4 \ge f_{\min}^2 \scm(\wh\beta_L)
r^2/4.
\end{eqnarray*}
Since the matrices $X^TX/n$ and $XX^T/n$ have the same maximal
eigenvalues,
$$
\frac{1}{n^2}({\bmf}-X\wh \beta_L)^T XX^T ({\bmf}-X\wh \beta_L)\le
\frac{\phi_{\max}}{n}|{\bmf}-X\wh \beta_L|_2^2 =
\phi_{\max}\|f-\wideth f_L\|_n^2
$$
and we deduce (\ref{Lassom}) from the last two displays.
\end{proof}
\begin{cor}\label{corx}
Let the assumptions of Lemma \ref{ll1} be satisfied and
$\|f_j\|_n=1, j=1,\dots,M$. Consider the linear regression model
$\bmy=X\bb+\bmw$. Then, with probability at least $1-M^{1-A^2/8}$,
we have
$$
|\bD_{J_0^c}|_1 \le 3|\bD_{J_0}|_1
$$
where $J_0=J(\beta)$ is the set of non-zero coefficients of $\beta$
and $\bD=\wh\beta_L-\beta$.
\end{cor}
\begin{proof} Use the first inequality in (\ref{onA}) and the fact
that $f=f_\beta$ for the linear regression model.
\end{proof}
\begin{lemma}\label{ll2}
Let $\beta\in\RR^M$ satisfy the Dantzig constraint
\begin{equation*}\begin{split}
\Big|\frac{1}{n}D^{-1/2}X^T(y-X\beta)\Big|_\infty \le r
\end{split}\end{equation*}
and set $\bD=\wh\beta_D-\beta$, $J_0=J(\beta)$. Then
\begin{equation}\begin{split}\label{ll22}
|\bD_{J_0^c}|_1\le |\bD_{J_0}|_1.
\end{split}\end{equation}
Further, let the assumptions of Lemma \ref{ll1} be satisfied with
$A>\sqrt{2}$. Then with probability of at least $1-M^{1-A^2/2}$ we
have
\begin{equation}\begin{split}\label{ll23}
\Big|\frac{1}{n}X^T({\bmf}-X\widehat\beta_D)\Big|_\infty \le
2rf_{\max}.
\end{split}\end{equation}
\end{lemma}

\begin{proof}[Proof of Lemma \ref{ll2}] Inequality(\ref{ll22}) follows immediately
from the definition of Dantzig selector, cf. \cite{ct07}. To prove
(\ref{ll23}) consider the event
\[ \B = \left\{\Big|\frac{1}{n}D^{-1/2}X^TW\Big|_\infty \le r\right\}
=\bigcap_{j=1}^M \left\{  | V_j | \leq r_{n,j} \right\}.\]
Analogously to (\ref{tail}), $\PP \{ \B^c\}  \le M^{1-A^2/2}$. On
the other hand, $\bmy={\bmf }+ \bmw$ and using the definition of Dantzig
selector it is easy to see that (\ref{ll23}) is satisfied on $\B$.

\end{proof}


\begin{proof}[Proof of Theorem \ref{th1}] Set $\bD=\wh\beta_L-\wh\beta_D$. We have
\begin{eqnarray*}
\frac{1}{n}|{\bmf}-X\wh \beta_L|_2^2 &= \frac{1}{n}|{\bf f}-X\wh
\beta_D|_2^2 -\frac{2}{n}\bD^TX^T({\bf f}-X\wh \beta_D) +
\frac{1}{n}|X \bD|^2_2.
\end{eqnarray*}
This and (\ref{ll23}) yield
\begin{equation}\label{e1}
\begin{split}
\| \wideth f_D -f\|_n^2 &\le \| \wideth{f}_L -f \|_n^2 +
2|\bD|_1 \, \Big|\frac{1}{n}X^T({\bf f}-X\wh \beta_D)\Big|_\infty
- \frac{1}{n}|X \bD|^2_2\\
 &\le \| \wideth{f}_L -f \|_n^2 + 4f_{\max}r|\bD|_1 -
\frac{1}{n}|X \bD|^2_2
\end{split}\end{equation}
where the last inequality holds with probability at least
$1-M^{1-A^2/2}$. Since the Lasso solution $\wh \beta_L$ satisfies
the Dantzig constraint, we can apply Lemma \ref{ll2} with $
\beta=\wh \beta_L$, which yields
\begin{equation}\begin{split}\label{e2}
|\bD_{J_{0}^c}|_1 \le |\bD_{J_0}|_1
\end{split}\end{equation}
with $J_0=J(\wh \beta_L)$. By Assumption \RE s1 we get
\begin{equation}\begin{split}\label{e3}
\frac{1}{\sqrt{n}}|X \bD|_2 \ge \kappa |\bD_{J_{0}}|_2
\end{split}\end{equation}
where $\kappa=\kappa(s,1)$. Using (\ref{e2}) and (\ref{e3}) we
obtain
\begin{equation}\begin{split}\label{e4}
|\bD|_1 \le 2|\bD_{J_{0}}|_1\le 2 \scm^{1/2}(\wh \beta_L)\
|\bD_{J_{0}}|_2\le \frac{2 \scm^{1/2}(\wh \beta_L)}{\kappa\sqrt{n}}|X
\bD|_2.
\end{split}\end{equation}
Finally, from (\ref{e1}) and (\ref{e4}) we get that, with
probability at least $1-M^{1-A^2/2}$,
\begin{equation}\label{e5a}
\begin{split}
\| \wideth f_D -f\|_n^2 &\le \| \wideth{f}_L -f \|_n^2 +
\frac{8 f_{\max} r \scm^{1/2}(\wh \beta_L)}{\kappa\sqrt{n}}|X \bD|_2
- \frac{1}{n}|X \bD|^2_2\\
&\le \| \wideth{f}_L -f \|_n^2 + \frac{16 f_{\max}^2 r^2 \scm(\wh
\beta_L)}{\kappa^2},
\end{split}\end{equation}
where the RHS follows  \eqref{dcf},  \eqref{ll23}, and another application of \eqref{e4}. This proves one side of the inequality.

To show the other side of the bound on the difference, we act as in (\ref{e1}), up to the inversion of
roles of $\wh\beta_L$ and $\wh\beta_D$, and we use (\ref{dcf}).
This yields that, with probability at least $1-M^{1-A^2/8}$,
\begin{equation}\label{e5}
\begin{split}
\| \wideth f_L -f\|_n^2 &\le \| \wideth{f}_D -f \|_n^2 +
2|\bD|_1 \, \Big|\frac{1}{n}X^T({\bf f}-X\wh \beta_L)\Big|_\infty
- \frac{1}{n}|X \bD|^2_2\\
 &\le \| \wideth{f}_D -f \|_n^2 + 3f_{\max}r|\bD|_1 -
\frac{1}{n}|X \bD|^2_2.
\end{split}\end{equation}
This is analogous to (\ref{e1}). Paralleling now the proof leading
to \eqref{e5a} we obtain
\begin{equation}\label{e5aaa}
\begin{split}
\| \wideth f_L -f\|_n^2 &\le \| \wideth{f}_D-f \|_n^2 + \frac{9
f_{\max}^2 r^2 \scm(\wh \beta_L)}{\kappa^2}\,.
\end{split}\end{equation}
The theorem now follows from (\ref{e5a}) and (\ref{e5aaa}).
\end{proof}



\begin{proof}[Proof of Theorem \ref{th2}] Set again $\bD=\wh\beta_L-\wh\beta_D$. We
apply (\ref{onA}) with $\beta=\wh\beta_D$ which yields that,
with probability at least $1-M^{1-A^2/8}$,
\begin{equation}\begin{split}\label{e6}
|\bD|_1 \le 4 |\bD_{J_0}|_1 + \| \wideth{f}_D -f \|_n^2 /r
\end{split}\end{equation}
where now $J_0=J(\wh \beta_D)$. Consider the two cases: (i) $\|
\wideth{f}_D -f \|_n^2
>2r|\bD_{J_0}|_1$ and (ii) $\| \wideth{f}_D -f \|_n^2
\le 2r|\bD_{J_0}|_1$. In case (i) inequality (\ref{e5}) with
$f_{\max}=1$ immediately implies
\begin{equation*}\begin{split}
\| \wideth{f}_L -f \|_n^2\le 10\| \wideth{f}_D -f \|_n^2
\end{split}\end{equation*}
and the theorem follows. In case (ii) we get from (\ref{e6}) that
\begin{equation*}\begin{split}
|\bD|_1 \le 6 |\bD_{J_0}|_1
\end{split}\end{equation*}
and thus $|\bD_{J_0^c}|_1 \le 5 |\bD_{J_0}|_1$. We can therefore
apply Assumption \RE s5 which yields, similarly to (\ref{e4}),
\begin{equation}\begin{split}\label{e8}
|\bD|_1 \le 6 \scm^{1/2}(\wh \beta_D)\ |\bD_{J_{0}}|_2\le \frac{6
\scm^{1/2}(\wh \beta_D)}{\kappa\sqrt{n}}|X \bD|_2
\end{split}\end{equation}
where $\kappa=\kappa(s,5)$. Plugging (\ref{e8}) into (\ref{e5}) we
finally get that, in case (ii),
\begin{equation}\begin{split}
\| \wideth f_L -f\|_n^2 &\le \| \wideth{f}_D -f \|_n^2 +
\frac{18\, r \scm^{1/2}(\wh \beta_D)}{\kappa\sqrt{n}}|X \bD|_2
- \frac{1}{n}|X \bD|^2_2\\
&\le\| \wideth{f}_D -f \|_n^2 +\frac{ 81\,r^2 \scm(\wh\beta_D)}
{\kappa^2}\,.
\end{split}
\end{equation}
\end{proof}


\begin{proof}[Proof of Theorem \ref{th-ora-las}] Fix an arbitrary $\beta\in\RR^M$ with
$\scm(\beta)\le s$. Set $\bD=D^{1/2}(\wh\beta_L-\beta)$,
$J_0=J(\beta)$. On the event $\A$, we get from the first inequality
in (\ref{onA}) that
\begin{equation}\label{onA1}
\begin{split}
\| \wideth f_L - f\|_n^2+r|\bD|_1
 &\le \| {\sf
f}_{\beta} - f\|_n^2 + 4r\sum_{j\in
J_0}\|f_j\|_{n} |\widehat \beta_{j,L} - \beta_j|\\
 & =\| {\sf f}_{\beta} - f\|_n^2 + 4r|\bD_{J_0}|_1,
 \end{split}
\end{equation}
and from the second inequality in (\ref{onA}) that
\begin{equation}\begin{split}\label{onA2}
\| \wideth f_L - f\|_n^2
&\le \| {\sf f}_{\beta} - f\|_n^2 + 4r
\sqrt{\scm(\beta)}\,|\bD_{J_0}|_2.
\end{split}\end{equation}
Consider separately the cases
\begin{equation}\begin{split}\label{onA3}
4r|\bD_{J_0}|_1  &\le \eps\|
 {\sf f}_{\beta} - f\|_n^2
\end{split}\end{equation}
and
\begin{equation}\begin{split}\label{onA4}
\eps \| {\sf f}_{\beta} - f\|_n^2 < 4r|\bD_{J_0}|_1.
\end{split}\end{equation}
In case (\ref{onA3}), the result of the theorem trivially follows
from (\ref{onA1}). So, we will only consider the case (\ref{onA4}).
All the subsequent inequalities are valid on the event $\A\cap \A_1$
where $\A_1$ is defined by (\ref{onA4}). On this event we get from
(\ref{onA1}) that
\begin{equation*}
\begin{split}
|\bD|_1 &\le 4(1+1/\eps)|\bD_{J_0}|_1
\end{split}
\end{equation*}
which implies $|\bD_{J_0^c}|_1 \le (3+4/\eps)|\bD_{J_0}|_1$. Then
for $\bD'=D^{-1/2}\bD=\wh\beta_L-\beta$ we have $|\bD_{J_0^c}'|_1
\le (3+4/\eps)(f_{\max}/f_{\min})|\bD_{J_0}'|_1$. We now use
Assumption \RE s{(3+4/\eps)f_{\max}/f_{\min}}. This yields
\begin{equation*}
\begin{split}
\kappa^2 |\bD_{J_0}'|_2^2 &\le
 \frac{1}{n}|X\bD'|_2^2 =
 \frac{1}{n}(\wh\beta_L-\beta)^T X^TX (\wh\beta_L-\beta)
 \\
 &= \| \wideth f_L -{\sf
 f}_{\beta}\|_n^2
 \end{split}
\end{equation*}
where $\kappa=\kappa(s,(3+4/\eps)f_{\max}/f_{\min})$. Combining this
with (\ref{onA2}) and using that $|\bD_{J_0}|_2 \le
f_{\max}|\bD_{J_0}'|_2$ we find
\begin{equation*}\label{onA7}
\begin{split}
\| \wideth f_L - f\|_n^2
 &\le \|
{\sf f}_{\beta} - f\|_n^2 +
4rf_{\max}\kappa^{-1}\sqrt{\scm(\beta)}\,\| \wideth f_L -{\sf
f}_{\beta}\|_n
\\
&\le \| {\sf f}_{\beta} - f\|_n^2 +
4rf_{\max}\kappa^{-1}\sqrt{\scm(\beta)}\left(\| \wideth f_L -f\|_n +
\|{\sf f}_{\beta}-f\|_n\right).
\end{split}
\end{equation*}
This inequality is of the same form as (A.4) in \cite{btw07}. A
standard decoupling argument as in \cite{btw07} using inequality
$2xy\le x^2/b + by^2$ with $b>1$, $x=r\kappa^{-1}\sqrt{\scm(\beta)}$,
and $y$ being either $\| \wideth f_L -f\|_n$ or $\|{\sf
f}_{\beta}-f\|_n$ yields that
\begin{eqnarray*}\label{5}
\| \wideth{f}_L -f\|_n ^2 \le  \frac{b+1}{b-1} \| {\sf f}_{\beta}
- f\|_n^2 + \frac{8b^2f_{\max}^2}{(b-1)\kappa^2} r^2 \scm(\beta), \ \
\forall \ b>1.
\end{eqnarray*}
Taking $b=1+2/\eps$ in the last display finishes
the proof of the theorem.\\
\end{proof}



\begin{proof}[Proof of Proposition \ref{th-ora-dan}] Due to the weak sparsity assumption there
exists $\bar \beta\in\RR^M$ with $\scm(\bar \beta)\le s$ such that
$\| {\sf f}_{\bar \beta} - f\|_n^2 \le {C_0
f_{\max}^2r^2}\kappa^{-2}\scm(\bar \beta)$ where $\kappa=\kappa(s,
3+4/\eps)$ is the same as in Theorem \ref{th-ora-las}. Using this
together with Theorem \ref{th-ora-las} and (\ref{Lassom}) we obtain
that, with probability at least $1-M^{1-A^2/8}$,
 \begin{eqnarray*}
\scm(\wh \beta_L) \le C_1(\eps)\scm(\bar \beta) \le C_1(\eps) s.
\end{eqnarray*}
This and Theorem \ref{th1} imply
 \begin{eqnarray*}
\| \wideth f_D -f\|_n^2 \le \| \wideth{f}_L -f\|_n ^2 + \,\frac{
16C_1(\eps) f_{\max}^2 A^2 \sigma^2}{\kappa_0^2}\, \left(\frac{s\log
M}{n}\right)\,
\end{eqnarray*}
where $\kappa_0= \kappa(\max (C_1(\eps),1)s, 3+4/\eps)$. Applying
Theorem \ref{th-ora-las} once again we get the result.
\end{proof}


\begin{proof}[Proof of Theorem \ref{th4}] Set $\bD=\widehat\bb_D-\bb^*$ and
$J_0=J(\bb^*)$. Using Lemma \ref{ll2} with $\beta=\bb^*$ we get
that on the event $\B$ (i.e., with probability at least
$1-M^{1-A^2/2}$): (i) $\frac{1}{n}|X^T X \bD|_\infty\le 2r$, and
(ii) inequality (\ref{l21}) holds with $c_0=1$. Therefore, on $\B$
we have
\begin{equation}\label{dl5}
\begin{split}
\frac{1}{n}|X \bD|_2^2&=\frac{1}{n}\, \bD^T X^T X \bD
\\
&\le \frac{1}{n}\, \Big|X^T X \bD\Big|_\infty\,|\bD|_1
\\
&\le 2r \Big(|\bD_{J_0}|_1 + |\bD_{J_0^c}|_1\Big)
\\
&\le 2(1+c_0)r |\bD_{J_0}|_1
\\
&\le 2(1+c_0)r\sqrt{s}\,|\bD_{J_0}|_2 = 4r \sqrt{s}\,|\bD_{J_0}|_2
\end{split}
\end{equation}
since $c_0=1$. From Assumption \RE s1 we get that
$$
\frac{1}{n}|X \bD|_2^2 \ge \kappa^2 |\bD_{J_0}|_2^2
$$
where $\kappa=\kappa(s,1)$. This and (\ref{dl5}) yield that,  on
$\B$,
\begin{equation}\begin{split}\label{dl6}
\frac{1}{n}|X \bD|_2^2 \le 16 r^2 s/\kappa^2, \quad |\bD_{J_{0}}|_2
\le 4 r \sqrt{s}/\kappa^2.
\end{split}\end{equation}
The first inequality in (\ref{dl6}) implies (\ref{th44}). Next,
(\ref{th43}) is straightforward in view of the second inequality in
(\ref{dl6}) and of the following relations (with $c_0=1$):
\begin{equation}\label{dl6a}
|\bD|_1= |\bD_{J_{0}}|_1 + |\bD_{J_{0}^c}|_1 \le
(1+c_0)|\bD_{J_{0}}|_1 \le (1+c_0)\sqrt{s}|\bD_{J_{0}}|_2
\end{equation}
that hold on $\B$. It remains to prove (\ref{th42}). It is easy to
see that the $k$th largest in absolute value element of
$\bD_{J_{0}^c}$ satisfies $|\bD_{J_{0}^c}|_{(k)}\le
|\bD_{J_{0}^c}|_1/k$. Thus
\begin{equation*}
|\bD_{J_{01}^c}|_2^2 \le |\bD_{J_{0}^c}|_1^2 \sum_{k\ge
m+1}\frac1{k^2} \le \frac1{m}|\bD_{J_{0}^c}|_1^2
\end{equation*}
and since (\ref{l21}) holds on $\B$ (with $c_0=1$) we find
$$
|\bD_{J_{01}^c}|_2 \le \frac{c_0|\bD_{J_{0}}|_1}{\sqrt{m}} \le
c_0|\bD_{J_{0}}|_2 \,\sqrt{\frac{s}{m}}\le
c_0|\bD_{J_{01}}|_2\,\sqrt{\frac{s}{m}}\ .
$$
Therefore, on $\B$,
\begin{equation}\label{chain0}
|\bD|_2\le \left(1+c_0\sqrt{\frac{s}{m}}\,\right)|\bD_{J_{01}}|_2.
\end{equation}
On the other hand, it follows from (\ref{dl5}) that
$$
\frac{1}{n}|X \bD|_2^2 \le 4 r \sqrt{s}\,|\bD_{J_{01}}|_2.
$$
Combining this inequality with Assumption \REm sm1 we obtain
that, on $\B$,
\begin{eqnarray*}
|\bD_{J_{01}}|_2 \le 4 r \sqrt{s}/\kappa^2.
\end{eqnarray*}
Recalling that $c_0=1$ and applying the last inequality together
with (\ref{chain0}) we get
\begin{equation}\label{chain2}
|\bD|_2^2\le 16 \left(1+c_0\sqrt{\frac{s}{m}}\,\right)^2 (r\sqrt{
s}/\kappa^2)^2.
\end{equation}
It remains to note that (\ref{th42}) is a direct consequence of
(\ref{th43}) and (\ref{chain2}). This follows from the fact that
inequalities $\sum_{j=1}^M a_j\le b_1$ and $\sum_{j=1}^M a_j^2\le
b_2$ with $a_j\ge0$ imply
$$
\sum_{j=1}^M a_j^p = \sum_{j=1}^M a_j^{2-p} a_j^{2p-2} \le
\left(\sum_{j=1}^M a_j\right)^{2-p}\left(\sum_{j=1}^M
a_j^2\right)^{p-1} \le b_1^{2-p} b_2^{p-1},\quad \forall \, 1<p\le
2.
$$

\end{proof}


\begin{proof}[Proof of Theorem \ref{th5}] Set $\bD=\widehat\bb_L-\bb^*$ and
$J_0=J(\bb^*)$. Using (\ref{onA}) where we put $\beta=\bb^*$,
$r_{n,j}\equiv r$ and $\|{\sf f}_\beta -f\|_n= 0$ we get that, on
the event $\A$,
\begin{equation}\begin{split}\label{dl8a}
\frac{1}{n}|X \bD|_2^2 \le 4r\sqrt{s}|\bD_{J_0}|_2
\end{split}\end{equation}
and (\ref{l21}) holds with $c_0=3$
on the same event. Thus, by Assumption \RE s3 and the last
inequality we obtain that, on $\A$,
\begin{equation}\begin{split}\label{dl8}  \frac{1}{n}|X \bD|_2^2 \le
16r^2 s/\kappa^2, \quad|\bD_{J_{0}}|_2\le 4r\sqrt{s}/\kappa^2
\end{split}\end{equation}
where $\kappa=\kappa(s,3)$. The first inequality here coincides with
(\ref{th54}). Next, (\ref{th55}) follows immediately from
(\ref{Lassom}) and (\ref{th54}). To show (\ref{th53}) it suffices to
note that on the event $\A$ the relations (\ref{dl6a}) hold with
$c_0=3$, to apply the second inequality in (\ref{dl8}) and to use
(\ref{tail}).

Finally, the proof of (\ref{th52}) follows exactly the same lines as
that of (\ref{th42}): the only difference is that one should set
$c_0=3$ in (\ref{chain0}), (\ref{chain2}), as well as in the display
preceding (\ref{chain0}).

\end{proof}

\end{document}